\documentclass{article}
\usepackage{amstext}
\usepackage{graphicx}
\usepackage{color}
\usepackage{amsmath}
\usepackage{amstext}
\usepackage{amsfonts}
\usepackage{amssymb}

\definecolor{Red}{cmyk}{0,1,1,0}

\definecolor{Blue}{cmyk}{1,1,0,0}

\definecolor{Pink}{cmyk}{0,1,0,0}

\definecolor{Green}{cmyk}{1,0,1,0.5}

\newcommand{\ba}{\begin{array}}
	\newcommand{\ea}{\end{array}}
\newcommand{\be}{\begin{equation}}
	\newcommand{\ee}{\end{equation}}
\newcommand{\ben}{\begin{enumerate}}
	\newcommand{\een}{\end{enumerate}}

\newcommand{\R}{\mathbb{R}}

\newtheorem{teo}{Theorem}[section]

\newtheorem{cor}{Corollary}[section]

\newtheorem{exer}{Exercise}[section]

\title{
	{\large{ \bf{The Diffusive Behavior of Solutions to the Linear Damped Wave
	Equation: an Undergraduate D.I.Y. Classnote}}}}

\date{}

\begin{document}
\maketitle

\centerline{\scshape Gast\~ao A. Braga}
\medskip
{\footnotesize
 \centerline{Departamento de Matem\'atica}
  \centerline{Universidade Federal de Minas Gerais}
   \centerline{Caixa Postal 1621, Belo Horizonte, 30161-970, Brazil}
} 

\medskip

\centerline{\scshape Ant\^onio Marcos da Silva}
\medskip
{\footnotesize
 \centerline{Departamento de Matem\'atica}
  \centerline{Universidade Federal de Ouro Preto}
   \centerline{R. Diogo de Vasconcelos, 122, Pilar, 35400-000, Brazil}
}

\medskip

\centerline{\scshape Jussara de Matos Moreira}
\medskip
{\footnotesize
\centerline{Departamento de Matem\'atica}
  \centerline{Universidade Federal de Minas Gerais}
   \centerline{Caixa Postal 1621, Belo Horizonte, 30161-970, Brazil}
}

\def\l{\lambda}

\baselineskip = 22pt

\maketitle

\begin{abstract}
Despite of the fact that the Damped Wave and the
Heat equations
describe phenomena of distinct nature, it is amazing that their solutions are
related in the limit as $t \to \infty$. The aim of this note is to explain to
undergraduate 
students, with a good calculus background,
how the relation between these solutions is established.
We follow a ``do it yourself'' strategy and the students
are invited to do the suggested exercises in order to
understand  the content of this note.
\end{abstract}
\clearpage

\section{\large{Introduction}}
\label{sec:intr}

Consider the following Partial Differential Equation (PDE)
\begin{equation}
\label{equa:tele-equa}
\mu u_t + u_{tt} - u_{xx} = 0,
\end{equation}
where $u=u(x,t)$, $x, t \in \R$ and $\mu\geq 0$. When $\mu=0$, this
is the one-dimensional wave equation, a classical example of a
{\em hyperbolic equation}.
Its solutions are the superposition of two travelling waves, one to the right and one to the left,
both propagating with velocity one, see \cite{bib:strauss}.
Jump discontinuities at
time zero will also propagate over characteristic curves with velocity one.
When $\mu>0$, the case we are interested in, the equation
is still hyperbolic and its solutions retain the properties of the $\mu=0$ solutions.
It is called 
the {\em Damped Wave Equation } (DWE),
or the {\em Telegraph Equation} (TE), see \cite{bib:cour-hilb} for instance.


\begin{exer}
\label{exer:ener-time-deri}
Multiplying both sides of (\ref{equa:tele-equa}) by $u_t$, integrating over $\R$ and assuming all operations
are allowed, conclude that
\begin{equation}
\label{equa:ener-time-deri}
\partial_t \left(\int_\R\frac{1}{2}[u_t^2 + u_x^2]dx \right)= -\mu \int_\R u_t^2 dx.
\end{equation}
\end{exer}
The integral on the left hand side (lhs) of (\ref{equa:ener-time-deri}) is the wave`s total energy.
If $\mu >0$ and if $u_t$ is not identically zero then the right hand side (rhs) of (\ref{equa:ener-time-deri}) is negative
implying that the wave's total energy decreases with time, not being conserved.
As we will see later on, the solutions to 
the DWE are also
a superposition of left and right travelling waves but, due to the damping term $\mu u_t$,
their amplitudes will diminish with time.

On the other hand, the {\em Heat} or {\em Diffusion Equation}, (HE) or (DE),
\begin{equation}
\label{equa:difu-equa}
\mu u_t  - u_{xx} = 0,
\end{equation}
where $u=u(x,t)$, $t>0, x \in \R$ and $\mu> 0$, is a classical example of a
{\em parabolic} PDE.
In \eqref{equa:difu-equa}, $\sigma = 1/\mu$ is the diffusion coefficient.
\begin{exer}
\label{exer:diff-equa-fund-solu}
For $t>0$ and $x\in\R$, show that
\begin{equation}
\label{equa:diff-equa-fund-solu}
K(x, t) = \frac{1}{\sqrt{t}}f^*_{\mu}\left(\frac{x}{\sqrt{t}}\right),
\end{equation}
where
\begin{equation}
\label{eq:dist-gaus}
f^*_{\mu}(x) = \sqrt{\frac{\mu}{4 \pi}}e^{- \mu \frac{x^2}{4}},
\end{equation}
is a solution to (\ref{equa:difu-equa}).
\end{exer}
We observe that $f^*_{\mu}(x)$, defined by (\ref{eq:dist-gaus}),
is the probability density function of a zero mean Gaussian random
variable with $\sqrt{2/\mu}$ variance.
\begin{exer}
\label{exer:mean-zero-gaus}
Show that
\begin{equation}
\label{equa:mean-zero-gaus}
\int_\R  f^*_{\mu}(x)dx = 1, \,\,\, \int_\R x f^*_{\mu}(x)dx = 0, \,\,\, \int_\R x^2 f^*_{\mu}(x)dx = \frac{2}{\mu}.
\end{equation}
\end{exer}

$K(x,t)$, given by \eqref{equa:diff-equa-fund-solu}, is said to be a {\em fundamental solution} to (\ref{equa:difu-equa}).
The following properties for $K(x,t)$ 
are easily verified:
\begin{exer}
\label{exer:fund-solu-prop}
For $t>0$ and $x\in\R$, verify that: $1)$  $K(x,t)$ is a $C^\infty$ function of $x$;
$2)$ $K(x,t)$ is scaling invariant, i.e.,
\be
\label{eq:fund-solu-resc}
\sqrt{t}K(\sqrt{t}x, t) = f^*_{\mu}(x).
\ee
\end{exer}
It turns out that any solution to the
Initial Value Problem (IVP), with $\mu >0$ and
a continuous $f(x)$,
\begin{equation}
\left\{
\begin{array}{ccccccc}
\mu u_t  - u_{xx} = 0, \ x \in \mathbb{R}, \ t > 0,\\
u(x, 0) = f(x),
\label{equa:ipv}
\end{array}
\right.
\end{equation}
retains the properties of the
$K(x, t)$, stated in Exercise \ref{exer:fund-solu-prop},
in the sense that: $1)$ for $t>0$, $u(\cdot,t)$ is a $C^\infty(\R)$ function
even if $u(x,t)$ has a jump discontinuity
at time $t=0$. In this case, we say that the
discontinuities at time $t=0$ are instantaneously smoothed out at any later time $t>0$;
$2)$ solutions to the IVP (\ref{equa:ipv}), as a function of $t$,
decay and spread out at rates $1/\sqrt t$ and $\sqrt t$, respectively,  i.e.,
identity \eqref{eq:fund-solu-resc} holds in the limit as $t\to\infty$,
see identity \eqref{eq:resc-solu-limi}.
The above two properties are a straightforward consequence of the well known
integral representation formula \cite{bib:strauss}.
\be
\label{eq:inte-repr}
u(x,t) = \int_\R K(x-y, t) f(y) dy
\ee
which holds for solutions to the IVP \eqref{equa:ipv} where $f(x)$
is a bounded and continuous function except possibly for
a finite number of jump discontinuities.
\begin{exer}
\label{exer:ivp-solu-prop}
Conclude, from (\ref{eq:inte-repr}), that the solution
$u(x,t)$ of the IVP \eqref{equa:ipv} with the above specified $f$,
is $C^{\infty}$ as a function of $x$.
\end{exer}
\begin{exer}
\label{exer:heat-equa-resc-solu-limi}
Conclude, from (\ref{eq:inte-repr}), that the solution
$u(x,t)$ of the IVP \eqref{equa:ipv} with the above specified $f$
such that $\int_\R |f(x)|dx <\infty$, satisfies
\be
\label{eq:resc-solu-limi}
 \lim_{t\to\infty}\sqrt{t}u(\sqrt{t}x, t) = M f^*_{\mu}(x),
\ee
where
\be
\label{eq:pre-fact}
M = \int_\R f(x)dx
\ee
and $f^*_{\mu}(x)$ is given by (\ref{eq:dist-gaus}).
 \end{exer}
The limit \eqref{eq:resc-solu-limi} expresses the fact that
identity \eqref{equa:diff-equa-fund-solu}, which is satisfied for the kernel $K(x,t)$,
holds asymptotically for the solution $u(x,t)$ of the IVP \eqref{equa:ipv} so that
$u(x,t)$ decays and spreads out at rates $1/\sqrt t$ and $\sqrt t$, respectively,
having  $f^*_{\mu}(x)$ as its {\em profile function}, in the limit $t\to\infty$.
$M$, given by \eqref{eq:pre-fact}, is the {\em prefactor}.
	
For large values of $t$, (\ref{eq:resc-solu-limi}) can be
rephrased as
\begin{equation}
\label{eq:comp-assi-difu-equa}
u(x, t) \approx \frac{M}{\sqrt{t}}f^*_{\mu}\left(\frac{x}{\sqrt{t}}\right).
\end{equation}
The above notation means that the two functions in (\ref{eq:comp-assi-difu-equa}) are asymptotically equivalent, that is, their ratio tends to one, when $t$ goes to infinity.
Despite of the fact that equations (\ref{equa:tele-equa}) and (\ref{equa:difu-equa})
describe phenomena of distinct nature, it is amazing that their solutions are
asymptotically related as $t \to \infty$. The aim of this note is to explain to
undergraduate mathematics, science and engineering students how the solutions to the above two equations
are connected.

Before presenting the main theorem and its proof, we provide a heuristic argument to highlight the intuition that distinct space and time scaling is behind the explanation for the changing mechanism from hyperbolic to parabolic behavior. To present this heuristic reasoning, let $u(x,t)$ be a solution to the following Cauchy Problem (CP)
\begin{equation}
\label{eq:cauc-prob}
\left\{
\begin{array}{ccccccc}
\mu u_t + u_{tt} - u_{xx} = 0, \ x \in \mathbb{R}, \ t > 0,\\
u(x, 0) = f(x), \\
u_t(x, 0) = g(x).
\end{array}
\right.
\end{equation}
Our purpose here is to prove that, if $u(x,t)$ is a solution to \eqref{eq:cauc-prob}, then the limit expressed in \eqref{eq:resc-solu-limi} holds, with the prefactor $M$ given by
 \begin{equation}
\label{eq:pre-fato}
 M = \int_{\mathbb{R}} \left[f(x) + \frac{1}{\mu}g(x) \right] dx.
 \end{equation}
Now define
\begin{equation}
\label{eq:rees-solu}
v(x,t) \equiv L^{\frac{1}{2}}u(L^{\frac{1}{2}}x, Lt),
\end{equation}
where $L>1$.
We say that the above $v(x,t)$ is a rescaling of
$u(x,t)$.
\begin{exer}
\label{exer:rees-equa}
Show that $v(x,t)$, defined by \eqref{eq:rees-solu}, solves the
equation
\begin{equation}
\mu v_t + \frac{1}{L}v_{tt} - v_{xx} = 0.
\label{eq:rees-equa}
\end{equation}
\end{exer}
Assuming that $|v_{tt}(x,t)|$ is uniformly bounded for $x\in\R$ and $t>0$, the second term on
the left hand side of \eqref{eq:rees-equa} will be small if we choose
$L$ large enough. Then, for large $L$, it is reasonable to drop this term off
thus generating the diffusion equation \eqref{equa:difu-equa}. More precisely, 
we conclude that, for large values of $t$,
\be
\label{eq:u-v}
u(x, t) \approx v(x,t),
\ee 
where
$v(x,t)$ is the solution to the IVP
\begin{equation}
\left\{
\begin{array}{ccccccc}
\mu u_t  - u_{xx} = 0, \ x \in \mathbb{R}, \ t > 0,\\
u(x, 0) = f(x) + \frac{1}{\mu}g(x).
\label{equa:modi-ipv}
\end{array}
\right.
\end{equation}

The approximation \eqref{eq:u-v} reflects the surprisingly fact that damped propagating
waves decay and spread out with rates $1/\sqrt t$ and $\sqrt t$, respectively.
This formal argument or, equivalently, the approximation
\eqref{eq:u-v}, is rigorously translated into the following theorem
which we will prove in the next sections.


\begin{teo}
\label{teo:asym-beha-tele-equa}
If $u(x,t)$ is the solution to the Cauchy problem \eqref{eq:cauc-prob}, with $\mu> 0$, $f \in C^2_0(\mathbb{R})$ and $g \in C^1_0(\mathbb{R})$, then
\be
\label{eq:resc-limi-tele-equa}
 \lim_{t\to\infty}\sqrt{t}u(\sqrt{t}x, t) = M f^*_{\mu}(x),
\ee
where $M$ and $f^*_{\mu}(x)$ are given by \eqref{eq:pre-fato} and \eqref{eq:dist-gaus}, respectively.
\end{teo}

{\bf Remark:} Theorem \ref{teo:asym-beha-tele-equa} is a simpler version of more general theorems
which require advanced mathematical methods to be proven. In the papers
\cite{bib:karch, bib:hsia-liu-01, bib:hsia-liu-02, bib:mats-01, bib:nishihara}, the reader
will be able to check in which directions the above theorem can be generalized.

\section{Integral representations of solutions to \eqref{eq:cauc-prob}}
\label{sec:inte-repr-bess-func}

The DW Equation \eqref{equa:tele-equa} is classified as hyperbolic and the pair of
straight lines
\begin{equation*}
\alpha = x + t \ \mbox{e} \ \beta = x - t.
\end{equation*}
form its characteristic curves,  (see \cite{Arthur-W}).
According to \cite{Arthur-W}, the solution to the CP \eqref{eq:cauc-prob}
can be expressed as follows
	\begin{eqnarray}
		u(x, t) =  \frac{e^{-\frac{\mu}{2}t}}{2}& \Big[ f(x+t) + f(x-t) 
		+ \displaystyle{\int_{x-t}^{x+t}}f(\alpha)\frac{d}{d t} I_0\left(\frac{\mu}{2}\sqrt{t^2 - (x - \alpha)^2} \right) d \alpha \nonumber  \\
		& \left. + \displaystyle{\int_{x-t}^{x+t}}\left(g(\alpha) + \frac{\mu}{2} f(\alpha)\right) I_0\left(\frac{\mu}{2}\sqrt{t^2 - (x - \alpha)^2} \right) d \alpha \right],
		\label{sol1}
	\end{eqnarray}
where $I_n(x)$, for $n = 0, 1, 2, \cdots$, represents the modified Bessel function of order $n$,
given by
 \begin{equation}
\label{eq:modi-bess-func}
I_n(x) = i^{-n} J_n(ix) = \sum_{j = 0}^{\infty} \frac{1}{j!(j+n)!}\left(\frac{x}{2}\right)^{2j+n},
 \end{equation}
where $J_n(x)$ is the Bessel function of order $n$,
 \begin{equation}
\label{eq:bess-func}
J_n(x) =  \sum_{j = 0}^{\infty} \frac{(-1)^j}{j!(j+n)!}\left(\frac{x}{2}\right)^{2j+n}.
 \end{equation}

If $\mu=0$ in \eqref{eq:cauc-prob} then we obtain the CP for the Wave Equation, whose solution
is given by the D'Alembert's formula 
\begin{equation}
	\begin{array}{rl}
		u(x, t) =  \frac{1}{2}[f(x+t) + f(x-t)]  + \frac{1}{2}\displaystyle{\int_{x-t}^{x+t}}g(\alpha)  d \alpha.
		\label{eq:solu-dala}
	\end{array}
\end{equation}

\begin{exer}
Show that D'Alembert's formula can be recovered from \eqref{sol1} when $\mu = 0$.
\end{exer}

The representation formula \eqref{sol1}, for $\mu>0$, inherits the left and right
propagating waves structure of \eqref{eq:solu-dala}:
at the point $P=(x_0, t_0)$, the solution is the superposition of two waves, both
travelling with velocity one along the characteristics and both being exponentially
damped if $\mu>0$.  Besides that, the value $u(x_0, t_0)$ depends uniquely on the
values of $f(\alpha)$ e $g(\alpha)$, for $\alpha \in \left[x_0 - t_0, x_0 + t_0\right]$,
the {\em domain of dependence}.

\begin{exer}
\label{exer2}
Using, for $n \in \mathbb{Z}$, that $\frac{d}{d \xi} \left( \xi^nI_n(\xi) \right) = \xi^nI_{n-1}(\xi)$ and
$I_{-n} (\xi) = I_{n} (\xi)$, show that $I_0'(\xi) = I_{-1}(\xi) = I_1(\xi)$,
where $I_0'(z)$ means $\frac{d}{dz}I_0(z)$ and conclude that the first integral on the right hand side of (\ref{sol1}) can be rewritten as
\begin{equation}
\displaystyle{\int_{x-t}^{x+t}}\frac{\mu}{2}f(\alpha)\frac{ t}{\sqrt{t^2 - (x - \alpha)^2} }
I_1\left(\frac{\mu}{2}\sqrt{t^2 - (x - \alpha)^2} \right) d \alpha.
\label{sol2}
\end{equation}
\end{exer}

Notice that the integration interval $|x-\alpha|\leq t$, also expressed
as  $t^2 - (x - \alpha)^2\geq 0$, leads that the Bessel functions $I_0$ in \eqref{sol1} and $I_1$ in \eqref{sol2} are in fact real numbers.

\section{Rescaling}
\label{sec:rees}

For fixed $x \in \mathbb{R}$, $t > 0$ and $\mu>0$, our aim is to verify that the rescaling \eqref{eq:rees-solu},
applied to the representation \eqref{sol1}, yields the limit \eqref{eq:resc-limi-tele-equa}.
For $L > 1$, define
\begin{equation}
\label{eq:xi}
\xi = \xi(\alpha; L, x)  \equiv \frac{\mu}{2}\sqrt{L^2 - (L^{\frac{1}{2}}x - \alpha)^2}.
\end{equation}

\begin{exer}
Use \eqref{sol1} and exercise \ref{exer2} to obtain a representation for the rescaling $v(x,t)$ given by \eqref{eq:rees-solu} and show that, at $t=1$, this representation is given by
$$
L^{\frac{1}{2}}u(L^{\frac{1}{2}}x, L) =  L^{\frac{1}{2}} \frac{e^{-\frac{\mu}{2}L}}{2}
\Big[f(L^{\frac{1}{2}}x+L) + f(L^{\frac{1}{2}}x-L) +
\int_{L^{\frac{1}{2}}x-L}^{L^{\frac{1}{2}}x+L} \frac{\mu^2 L}{4\xi}f(\alpha)I_{1}(\xi) d \alpha
$$
\begin{equation}
 + \left. \displaystyle{\int_{L^{\frac{1}{2}}x-L}^{L^{\frac{1}{2}}x+L}}
\left(g(\alpha) + \frac{\mu}{2} f(\alpha)\right) I_0\left(\xi \right) d \alpha \right].
\label{solreesc00}
\end{equation}
\end{exer}

Notice that we can replace $L$ by $L/t$ in the rescaled function $v(x,t)$,
so that the results will hold for $L^{\frac{1}{2}}u(L^{\frac{1}{2}}y, L)$,
being enough to replace $y$ by $x/\sqrt t$ and multiply $u$ by
$1/\sqrt t$. That is why, from now on we will consider $t=1$ and therefore, we will use \eqref{solreesc00}
to analyze the behavior of  $L^{\frac{1}{2}}u(L^{\frac{1}{2}}x, L)$ when $L\gg 1$.
Furthermore, notice that the Bessel functions $I_0$ and $I_1$ in equation \eqref{solreesc00} are also real numbers.



Since $f(x)$ is a compact support function, then, $f(x)=0$ if $x\not\in I_f$, where
$I_f$ is its support. In particular, there exists $L_{0} = L_{0}(x)>1$ such that if $L>L_{0}$ then   $(L^{\frac{1}{2}}x\pm L)\not\in I_{f}$,
i.e., $f(L^{\frac{1}{2}}x+L) = 0=  f(L^{\frac{1}{2}}x-L)$. Therefore, if $L>L_0$, then, the right hand side of
\eqref{solreesc00} can be rewritten as
\begin{equation}
 \frac{\sqrt{L}e^{-\frac{\mu}{2}L}}{2}
\left[
\displaystyle\int_{\sqrt{L}x-L}^{\sqrt{L}x+L}\frac{\mu^2 L}{4\xi}f(\alpha)I_{1}(\xi) d \alpha  +
\displaystyle\int_{\sqrt{L}x-L}^{\sqrt{L}x+L}
\left[g(\alpha) + \frac{\mu}{2} f(\alpha)\right] I_0\left(\xi \right) d \alpha \right].
\label{solreesc}
\end{equation}

\section{Approximations for Bessel Functions}
\label{sec:apro-func-bess}

It follows from Theorem \ref{teo:assi-limi} that,
given $n \in \mathbb{Z}$, there exist positive constants $C$ and $\xi_0$ such that,
for all $\xi>\xi_0$,
\begin{equation}
\left|I_n(\xi) - \frac{1}{\sqrt{2\pi}}\frac{e^{\xi}}{ \sqrt \xi}\right|
\leq \frac{C}{\sqrt{ 2\pi}} \frac{e^{\xi}}{\xi}.
\label{assI}
\end{equation}

We want to use \eqref{assI} to estimate the integrals on \eqref{solreesc}. In order to do that, we must ensure that $\xi$, defined by (\ref{eq:xi}) and which appears as the argument of  $I_0$ and $I_1$,
will satisfy the condition $\xi>\xi_0$.
\begin{exer}
\label{exer:51}
Define $L_1=2\xi_0/\mu$ and $\overline{L} = \sqrt{L^2-L_1^2}$. Show that $\xi>\xi_0$ if and only if $L>L_1$ and $\alpha \in [L^{1/2}x-\overline{L}, L^{1/2}x+\overline{L}]$.
\end{exer}


\begin{exer}
\label{exer:52}
Show that, for $L>L_1$, there exists $L_2(x)$ such that, if $L>L_2$, then $f$ and $g$ vanish outside the interval $[L^{1/2}x-\overline{L}, L^{1/2}x+\overline{L}]$.
\end{exer}

From exercises \ref{exer:51} and \ref{exer:52}, if we take $L> \max\{L_0, L_1, L_2\}$ then
we are allowed to replace \eqref{assI} on \eqref{solreesc}. In what follows, we prove the limit (\ref{eq:resc-limi-tele-equa}) by showing that, as  $L\to\infty$, the term
$$
\left|
\frac{\sqrt{L}e^{-\frac{\mu}{2}L}}{2}
\left\{
\displaystyle{\frac{C}{\sqrt{2\pi}} \int_{L^{\frac{1}{2}}x-\overline{L}}^{L^{\frac{1}{2}}x+\overline{L}}}\frac{\mu^2 Le^{\xi}}{4\xi^2}f(\alpha)
d \alpha
+
\displaystyle{\frac{C}{\sqrt{2\pi}} \int_{L^{\frac{1}{2}}x-\overline{L}}^{L^{\frac{1}{2}}x+\overline{L}}}
\left[g(\alpha) + \frac{\mu}{2} f(\alpha)\right] \frac{e^{\xi}}{\xi} d \alpha \right\}
\right|
$$
goes to zero, while the term
$$
\frac{\sqrt{L}e^{-\frac{\mu}{2}L}}{2}
\left\{
\displaystyle{
\frac{1}{\sqrt{2\pi}} \int_{L^{\frac{1}{2}}x-\overline{L}}^{L^{\frac{1}{2}}x+\overline{L}}}\frac{\mu^2 Le^{\xi}}{4\xi^{3/2}}f(\alpha)
 d \alpha
+  \displaystyle{\frac{1}{\sqrt{2\pi}} \int_{L^{\frac{1}{2}}x-\overline{L}}^{L^{\frac{1}{2}}x+\overline{L}}}
\left[g(\alpha) + \frac{\mu}{2} f(\alpha)\right] \frac{e^{\xi}}{ \sqrt \xi} d \alpha \right\}
$$
goes to
$Mf^*_{\mu}\left(x\right)$,
where $f_{\mu}^*$ is the Gaussian distribution \eqref{eq:dist-gaus}
and $M$ is given by \eqref{eq:pre-fato}.

Defining
$$
T_\beta = \frac{\mu}{4\sqrt{2\pi}} \int_{L^{\frac{1}{2}}x-\overline{L}}^{L^{\frac{1}{2}}x+\overline{L}}
\left(L^{3/2}e^{-\frac{\mu}{2}L}\frac{e^{\xi}}{\xi^\beta}\right)\frac{\mu}{2} f(\alpha)d\alpha
$$
and
$$
S_\beta = \frac{1}{2\sqrt{2\pi}}\int_{L^{\frac{1}{2}}x-\overline{L}}^{L^{\frac{1}{2}}x+\overline{L}}
\left(L^{1/2}e^{-\frac{\mu}{2}L}\frac{e^{\xi}}{\xi^\beta}\right)\left[g(\alpha) + \frac{\mu}{2} f(\alpha)\right]d\alpha,
$$
we must prove that $|C(T_2+S_1)| \to 0$ and $T_{3/2} + S_{1/2} \to Mf_\mu^*$ as $L\to \infty$. Notice that, if we rewrite $\xi$ as
$$
\xi = \frac{\mu L}{2}\left[\sqrt{1 - \left(\frac{x}{L^{\frac{1}{2}}} - \frac{\alpha}{L}\right)^2} - 1\right]  +\frac{\mu L}{2},
$$
then, the $L-$dependence of the preceding integrands can be expressed in the form
\begin{equation}
\label{eq:L-depe-inte-01}
\left(\frac{ 2}{\mu}\right)^{\beta}L^{\gamma-\beta}
\left\{
 \frac{e^{\frac{\mu L}{2}\left[\sqrt{1 - \left(\frac{x}{L^{\frac{1}{2}}} -
\frac{\alpha}{L}\right)^2} - 1\right]}}{\left[1 - \left(\frac{x}{L^{\frac{1}{2}}} - \frac{\alpha}{L} \right)^2 \right]^{\beta/2}}
\right\},
\end{equation}
where $\gamma=3/2$ for $T_\beta$, $\gamma=1/2$ for $S_\beta$ and $\beta\in\{1/2,1,3/2,2\}$.

\begin{exer}
\label{exer:tayl-expa-squa-root}
Using the Taylor expansion with integral reminder, show that
$$
\left|\sqrt{1-x} - 1 + \frac{1}{2}x\right| \leq \max_{[-x,0]}\left\{\frac{1}{4}(1+t)^{-3/2}\right\}x^2,\,\,\, |x|<1.
$$
\end{exer}
\begin{exer}
\label{exer:squa-root-limi}
Using Exercice \ref{exer:tayl-expa-squa-root}, show that
$$
\lim_{L\to\infty} \frac{\mu L}{2}\left[\sqrt{1 - \left(\frac{x}{L^{\frac{1}{2}}} -
\frac{\alpha}{L}\right)^2} - 1\right] = -\mu\frac{x^2}{4}.
$$
\end{exer}
From Exercise \ref{exer:squa-root-limi}, it follows that the term enclosed in braces in \eqref{eq:L-depe-inte-01} converges to $e^{-\mu\frac{x^2}{4}}$ as $L\to\infty$. Applying the Dominated Convergence Theorem \cite{bib:rudin}, one readily verifies that
 $|C(T_2+S_1)| \to 0$, since the factor $L^{\gamma-\beta}$ reduces to $L^{3/2-2}$ for $T_2$ and $L^{1/2-1}$ for $S_1$, both of which vanish as $L\to\infty$. Moreover, for $T_{3/2}$ and $S_{1/2}$, the factor simplifies to $L^{\gamma-\beta}=L^0=1$, so that, in the limit $L\to\infty$, $T_{3/2}+S_{1/2}$ converges to
$$
 \frac{1}{2\sqrt{2\pi}}\left(\frac{2}{\mu}\right)^{1/2}
e^{-\mu\frac{x^2}{4}}
\left\{ \int_{-\infty}^{\infty}\frac{\mu}{2}f(\alpha) d \alpha \right.
\left. + \int_{-\infty}^{\infty}
\left[g(\alpha) + \frac{\mu}{2} f(\alpha)\right]  d \alpha \right\} =
$$
$$
\sqrt{\frac{ \mu}{4\pi }}
e^{-\mu\frac{x^2}{4}}
\int_{-\infty}^{\infty}
\left[\frac{1}{\mu}g(\alpha) + f(\alpha)\right]  d \alpha = Mf_\mu^*(x).
$$

\appendix

\section{Bessel Functions}

From the series representation (\ref{eq:modi-bess-func}), one sees
that $I_n(x)$, $x\in \R$, is an even or odd function depending if
$n$ is even or odd, respectively.
Also, for fixed $(x,t)$ and for $\alpha \in [x-t, x+t]$, the argument of $I_0(x)$ and $I_1(x)$,
in the integral representation (\ref{sol1}), is non-negative, see Remark 2 of Section \ref{sec:inte-repr-bess-func}.
Said that, in the sequel we will prove the asymptotic behaviour (\ref{assI}):
\begin{teo}
\label{teo:assi-limi}
Consider the modified Bessel function $I_\nu(x)$, with $\nu\geq 0$ and $x > 0$.
Then
\begin{equation}
\label{eq:assi-limi}
\left|
\frac{\sqrt{2 \pi x}}{e^x} I_\nu(x) - 1		
\right|
\leq C(x),
\end{equation}
where
$$
C(x)\equiv \sqrt{\frac{\pi^3 x}{2}}e^{-x(1-\cos \delta)} +
 \sqrt{\frac{\pi^3}{2x}} (e^{-x}-e^{-2x}) +
e^{-2x}+\frac{2|\nu|}{ex\sigma}+\frac{128\sqrt{2}}{xe^2\sigma^{5/2}}+\sqrt{2}e^{- \frac{\delta^2 x}{4}},
$$
with $\sigma = 1- 4\delta^2>0$ and $0<\delta<1/2$.
\end{teo}
\begin{cor}
\label{cor:ratio-limi}
\begin{equation}
\label{eq:ratio-limi}
\lim_{x\to\infty}
\frac{\sqrt{2 \pi x}}{e^x} I_\nu(x) = 1.
\end{equation}
\end{cor}
{\bf Remark:} Corollary \ref{cor:ratio-limi} says that $I_\nu(x)$
and $\frac{e^x}{\sqrt{2 \pi x}}$ are asymptotically equivalent as
$x\to\infty$ while Theorem \ref{teo:assi-limi} says that the
difference $I_\nu(x) - \frac{e^x}{\sqrt{2 \pi x}}$
is a little order of $\frac{e^x}{\sqrt{2 \pi x}}$ as
$x\to\infty$. Therefore, inequality \eqref{eq:assi-limi} is a stronger
statement since it implies the limit \eqref{eq:ratio-limi}. 

{\bf Proof of Theorem \ref{teo:assi-limi}:} We make it explicitly the necessary estimates needed to
prove estimate (\ref{eq:assi-limi}.
The starting point is the following integral representation
\begin{equation}
\label{eq:mbf-inte-repr}
I_\nu(x) = \frac{1}{\pi}\int_0^\pi e^{x\cos t} \cos(\nu t) dt - \frac{\sin (\nu\pi)}{\pi}\int_0^\infty e^{-x\cosh t - \nu t} dt,
\end{equation}
see \cite{bib:watson,bib:mainardi}.

\begin{exer}
\label{exer:ap1}
Show that, for $a>0$, $\int_{-\infty}^\infty e^{-at^2} dt = \sqrt{\frac{\pi}{a}}$, then, use this to prove that, for $\nu \geq 0$,
$$
\int_0^\infty e^{-x\cosh t - \nu t} dt \leq \sqrt{\frac{\pi}{2 x}} e^{-x}.
$$
\end{exer}
From \eqref{eq:mbf-inte-repr} and exercise \ref{exer:ap1}, we get
\begin{equation}
\label{eq:uppe-boun-01}
\left|
I_\nu(x) - \frac{1}{\pi}\int_0^\pi e^{x\cos t} \cos(\nu t) dt
\right|
\leq
\sqrt{\frac{1}{2\pi x}} e^{-x}, \,\,\ x>0.
\end{equation}
Now we rewrite the integral on the lhs of  \eqref{eq:uppe-boun-01}
as a sum of two integrals, the second of which
being
$$
\int_{\pi/2}^{\pi} e^{x\cos t} \cos(\nu t) dt
=
\int_0^{\pi/2} e^{-x\sin u} \cos(\nu (u+\pi/2)) du.
$$
\begin{exer}
\label{exer:ap2}
Show that, if $u\in [0, \pi/2]$, then $u\geq \sin u \geq (2u)/\pi$ and use this to prove that
$$
\left| \int_{\pi/2}^{\pi} e^{x\cos t} \cos(\nu t) dt
\right| \leq
\pi \frac{1-e^{-x}}{2x}, \,\,\ x>0.
$$
\end{exer}
From (\ref{eq:uppe-boun-01}) and exercise \ref{exer:ap2}, we get
\begin{equation}
\label{eq:uppe-boun-02}
\left|
I_\nu(x) - \frac{1}{\pi}\int_0^{\pi/2} e^{x\cos t} \cos(\nu t) dt
\right|
\leq
\pi \frac{1-e^{-x}}{2x} +
\sqrt{\frac{1}{2\pi x}} e^{-x}, \,\,\ x>0.
\end{equation}
To handle the integral on the left-hand side of \eqref{eq:uppe-boun-02}, we decompose it into two parts:
\begin{equation}
\label{eq:uppe-boun-02-1}
\frac{1}{\pi}\int_0^{\pi/2} e^{x\cos t} \cos(\nu t) dt
=\frac{1}{\pi}\int_0^{\delta} e^{x\cos t} \cos(\nu t) dt
+
\frac{1}{\pi}\int_\delta^{\pi/2} e^{x\cos t} \cos(\nu t) dt,
\end{equation}
where $0< \delta < \pi/2$. We then choose $\delta$ sufficiently small so that, for $t\in[0,\delta]$, the approximation $x\cos t \approx x\left(1 - \tfrac{t^2}{2}\right)$ is valid, which will allow us to estimate
the first integral in the rhs of \eqref{eq:uppe-boun-02-1}.

\begin{exer}
\label{exer:ap3}
Show that the second integral in the rhs of \eqref{eq:uppe-boun-02-1} can be bounded above by $\pi e^{x\cos \delta}/2$.
\end{exer}
Later on, we will choose $\delta$ as a function of $x$ so that the appropriate
limits can be taken as $x\to\infty$. Having this in mind, we use exercise \ref{exer:ap3} to replace \eqref{eq:uppe-boun-02} by
\begin{equation}
\label{eq:uppe-boun-03}
\left|
I_\nu(x) - \frac{1}{\pi}\int_0^{\delta} e^{x\cos t} \cos(\nu t) dt
\right|
\leq \frac{\pi}{2} e^{x\cos \delta} +
\pi \frac{1-e^{-x}}{2x} +
\sqrt{\frac{1}{2\pi x}} e^{-x}, \,\,\ x>0.
\end{equation}

\begin{exer}
\label{exer:ap4}
Let $t\in [0, \delta]$, $0<\delta < 1$ and $R(t) \equiv \cos  t - (1 - t^2/2)$. Show that
$$
|R(t)| \leq (\cosh{1}) t^4.
$$
\end{exer}
Replacing $\cos t$ by $1 - t^2/2 + R(t)$ in \eqref{eq:uppe-boun-03} and multiplying both sides by $(\sqrt{2 \pi x})e^{-x}$, we get
$$
\left|
\frac{\sqrt{2 \pi x}}{e^x} I_\nu(x) -
\sqrt{\frac{2}{\pi}}\int_0^{\delta \sqrt{x}} e^{- \frac{u^2}{2} + xR\left(\frac{u}{\sqrt{x}}\right)} \cos\left(\nu \frac{u}{\sqrt{x}}\right) du
\right|
\leq
$$
\begin{equation}
\label{eq:uppe-boun-05}
\sqrt{\frac{\pi^3 x}{2}}e^{-x(1-\cos \delta)} +
 \sqrt{\frac{\pi^3}{2x}} (e^{-x}-e^{-2x}) +
e^{-2x}, \,\,\ x>0.
\end{equation}
Taking into account \eqref{eq:uppe-boun-05} and making use of the triangle inequality, to obtain the desired bound for $\left| \sqrt{2 \pi x}e^{-x} I_\nu(x) - 1 \right|$ in Theorem \ref{teo:assi-limi}, it remains to estimate the term
$$
\left|
\sqrt{\frac{2}{\pi}}\int_0^{\delta \sqrt{x}} e^{- \frac{u^2}{2} + xR\left(\frac{u}{\sqrt{x}}\right)} \cos\left(\nu \frac{u}{\sqrt{x}}\right) du
-1\right|= \sqrt{\frac{2}{\pi}}\left|V_x-\sqrt{\frac{\pi}{2}}\right|,
$$
where
$V_x \equiv \int_0^{\delta \sqrt{x}} e^{- \frac{u^2}{2} + xR\left(\frac{u}{\sqrt{x}}\right)} \cos\left(\nu \frac{u}{\sqrt{x}}\right) du$.
Our goal from now on is to establish an upper bound for the term $|V_x-\sqrt{\pi/2}|$, which we will accomplish by bounding
\begin{equation}
\label{bound1}
\left|V_x -\int_0^{\delta \sqrt{x}} e^{- \frac{u^2}{2} + xR\left(\frac{u}{\sqrt{x}}\right)} du\right|+
\left|\int_0^{\delta \sqrt{x}} e^{- \frac{u^2}{2} + xR\left(\frac{u}{\sqrt{x}}\right)} du- \sqrt{\frac{\pi}{2}}\right|.
\end{equation}
For this purpose, we shall employ the following two exercises.
\begin{exer}
\label{exer:ap5}
Show that, for $t \geq 0$, we have $|\cos t - 1|\leq t$.
\end{exer}

\begin{exer}
\label{exer:ap6}
Show that, for $t\in[0,\delta]$ and $0<\delta <1$, the inequality in exercise \ref{exer:ap4} can be improved to
$0\leq R(t) \leq  (\cosh{1})t^4$.

\end{exer}
From exercise \ref{exer:ap6} we get that $R(t)\leq 2t^4 \leq 2t^2\delta^2$. Therefore,
$$
- \frac{u^2}{2} \leq - \frac{u^2}{2} + xR\left(\frac{u}{\sqrt{x}}\right)
\leq -\left(1- 4\delta^2\right)\frac{u^2}{2}
$$
and combining this result with Exercise \ref{exer:ap5}, we can bound the first term in \eqref{bound1} as follows:
\begin{equation}
\label{bound2}
\left|\int_0^{\delta \sqrt{x}} e^{- \frac{u^2}{2} + xR\left(\frac{u}{\sqrt{x}}\right)}
\left[\cos\left(\nu \frac{u}{\sqrt{x}}\right)-1 \right]du\right|\leq
\frac{|\nu|}{\sqrt{x}} \int_0^{\infty} u e^{-\left(1- 4\delta^2\right)\frac{u^2}{2}} du.
\end{equation}
Before proceeding to the analysis of the above integral, we notice that, since
$$
\sqrt{\frac{\pi}{2}} =
\int_0^{\infty} e^{- \frac{u^2}{2}} du
=
\int_0^{\delta \sqrt{x}} e^{- \frac{u^2}{2}} du
+
\int_{\delta \sqrt{x}}^{\infty} e^{- \frac{u^2}{2}} du,
$$
we can rewrite the second term in \eqref{bound1} as
$$
\left|\int_0^{\delta \sqrt{x}} e^{- \frac{u^2}{2}}\left[e^{xR\left(\frac{u}{\sqrt{x}}\right)}-1\right] du-
\int_{\delta \sqrt{x}}^{\infty} e^{- \frac{u^2}{2}} du \right|.
$$

\begin{exer}
\label{exer:ap7}
Show that $e^t - 1\leq te^t$, for $t \geq 0$.
\end{exer}
Using exercises \ref{exer:ap6} and \ref{exer:ap7}
and since
$
0\leq
\int_{\delta \sqrt{x}}^{\infty} e^{- \frac{u^2}{2}} du
\leq
e^{- \frac{\delta^2 x}{4}}
\int_{0}^{\infty} e^{- \frac{u^2}{4}} du
=
\sqrt{\pi}e^{- \frac{\delta^2 x}{4}}
$, 
we obtain the following upper bound for the second term in \eqref{bound1}:
\begin{equation}
\label{bound3}
\frac{2}{x}\int_0^{\delta \sqrt{x}}u^4 e^{-(1- 4\delta^2)\frac{u^2}{2}}  du+\sqrt{\pi}e^{- \frac{\delta^2 x}{4}}.
\end{equation}
\begin{exer}
\label{exer:ap8}
Show that, for fixed positive constants $\alpha, \sigma$,
$$
\int_0^{\infty} u^{\alpha} e^{-\sigma\frac{u^2}{2}} du\leq \left(\frac{2\alpha}{\sigma e}\right)^{\alpha/2} \sqrt{\frac{\pi}{\sigma}}.
$$
\end{exer}
Notice that we can use exercise \ref{exer:ap8} to bound the integrals in \eqref{bound2} and \eqref{bound3}, with $\alpha \in \{1,4\}$, as long as $\delta$ is small enough so that $\sigma = 1- 4\delta^2>0$. Accordingly, our estimate for \eqref{bound1} takes the form
$$
\frac{\sqrt{2\pi}}{ex}\frac{|\nu|}{\sigma}+\frac{128}{xe^2}\frac{\sqrt{\pi}}{\sigma^{5/2}}+\sqrt{\pi}e^{- \frac{\delta^2 x}{4}}.
$$
Using the above estimate and \eqref{eq:uppe-boun-05}, we can finally obtain as a bound for $\left| \sqrt{2 \pi x}e^{-x} I_\nu(x) - 1 \right|$:
$$
\sqrt{\frac{\pi^3 x}{2}}e^{-x(1-\cos \delta)} +
 \sqrt{\frac{\pi^3}{2x}} (e^{-x}-e^{-2x}) +
e^{-2x}+\frac{2|\nu|}{ex\sigma}+\frac{128\sqrt{2}}{xe^2\sigma^{5/2}}+\sqrt{2}e^{- \frac{\delta^2 x}{4}},
$$
where $\sigma = 1- 4\delta^2>0$.


\end{document}